\documentclass{jnmp}

\usepackage{amsmath}

\JNMPnumberwithin{equation}{section}

\newtheorem{theorem}{Theorem}
\newtheorem{lemma}{Lemma}

\theoremstyle{definition}
\newtheorem{remark}{Remark}

\begin{document}

\renewcommand{\evenhead}{P~Redou}
\renewcommand{\oddhead}{Representations of the Conformal Lie Algebra}

\thispagestyle{empty}

\FirstPageHead{10}{2}{2003}{\pageref{redou-firstpage}--\pageref{redou-lastpage}}{Letter}

\copyrightnote{2003}{P~Redou}

\Name{Representations of the Conformal Lie Algebra \\
in the Space of Tensor Densities on the Sphere}
\label{redou-firstpage}

\Author{Pascal REDOU}

\Address{Institut Girard Desargues, Universit\'e Claude Bernard Lyon 1, \\
B\^atiment Braconnier (ex-101), 21 Avenue
Claude Bernard, \\
69622 Villeurbanne Cedex, France\\
E-mail: redou@enib.fr}

\Date{Received June 11, 2002;
Accepted August 2, 2002}

\begin{abstract}
\noindent
Let ${\mathcal F}_\lambda(\mathbb{S}^n)$ be the space of tensor densities on
$\mathbb{S}^n$ of degree $\lambda$. We consider this space as an induced module
of the nonunitary spherical series of the group
$\mathrm{SO}_0(n+1,1)$ and classify $(\mathrm{so}(n+1,1),\mathrm{SO}(n+1))$-simple and
unitary submodules of ${\mathcal F}_\lambda(\mathbb{S}^n)$
as a~function of $\lambda$.
\end{abstract}

\section{Introduction and main result}

Let ${\mathcal F}_\lambda(\mathbb{S}^n)$ be the space of tensor densities of degree
$\lambda\in{\mathbb C}$ on the sphere $\mathbb{S}^n$, that is, of smooth sections of the
line bundle
\[
\Delta_\lambda(\mathbb{S}^n)=
|\Lambda^nT^*\mathbb{S}^n|^{\otimes\lambda}
\]
on $\mathbb{S}^n$.
This space
plays an important r\^ole in geometric quantization
and, more recently, it has also been used in equivariant quantization
(see~\cite{DLO}). This space is endowed with a structure of
$\mathrm{Diff}(\mathbb{S}^n)$- and $\mathrm{Vect}(\mathbb{S}^n)$-module in the following way.
As a vector space, it is isomorphic to the space
${\mathcal C}^\infty_{\mathbb C}(\mathbb{S}^n)$ of smooth complex-valued functions;
the action of a~vector field
\[
Y=\sum_{i=1}^nY_i\frac{\partial}{\partial{}x_i}
\]
is given by the Lie derivative of degree $\lambda$
\begin{equation}
L_Y^\lambda(\varphi(x_1,\ldots,x_n))=
\sum_{i=1}^n\left(Y_i\frac{\partial{}
\varphi}{\partial{}x_i}+\lambda\,\frac{\partial
{}Y_i}{\partial{}x_i}\,\varphi\right)(x_1,\dots,x_n)
\label{LieDer}
\end{equation}
in any coordinate system.

The Lie algebra $\mathrm{so}(n+1,1)\subset\mathrm{Vect}(\mathbb{S}^n)$ of infinitesimal
conformal transformations, that we call the conformal Lie algebra, is
generated by the vector fields
\begin{gather}
X_i  = \frac{\partial}{\partial s_i}, \qquad
X_{ij}=
 s_i\frac{\partial}{\partial s_j}-
s_j\frac{\partial}{\partial s_i},\nonumber\\
X_0 =\sum_i s_i\frac{\partial}{\partial s_i},\qquad
\bar X_i
=\sum_j\left( s_j^2\frac{\partial}{\partial s_i}-
2s_is_j\frac{\partial}{\partial s_j}\right),
\label{confGenerators}
\end{gather}
where $(s_1,\ldots,s_n)$ are stereographic coordinates on the sphere
$\mathbb{S}^n$.

The space ${\mathcal F}_\lambda(\mathbb{S}^n)$ is naturally an $\mathrm{so}(n+1,1)$-module;
furthermore, the restriction of the action of the group
$\mathrm{Diff}(\mathbb{S}^n)$, defines the action of the subgroup $\mathrm{SO}(n+1)$ given
by the formula
\[
(k_0.f)(k)=f(k_0^{-1}k),\qquad {\rm where}\quad k_0\in
K, \ \ k\in\mathbb{S}^n\simeq\mathrm{SO}(n+1)/\mathrm{SO}(n).
\]
Therefore, ${\mathcal F}_\lambda(\mathbb{S}^n)$ is also a $\mathrm{SO}(n+1)$-module.

Given a Lie group $G$ and a compact subgroup $K\subset{G}$, let
$\mathfrak{g}$ and $\mathfrak{k}$ be the corresponding Lie algebras.
One calls
$(\mathfrak{g},K)$-module a complex vector space $E$ endowed with
actions of~$\mathfrak{g}$ and $K$ such that
\begin{enumerate}
\item
$({\rm Ad}k\cdot X)\cdot e=k\cdot X\cdot k^{-1}\cdot e\qquad \forall\;
k\in K, \quad X\in\mathfrak{g}, \quad e\in E$
\item
For all $e\in E$, the space $K\cdot e$ is finite-dimensional
(i.e., $e$ is a $K$-finite vector), the representation of $K$
in $F$ is continuous and one has for $X\in\mathfrak{k}$:
\[
X\cdot e=\frac{d}{dt}(\exp tX)\cdot e|_{t=0}.
\]
\end{enumerate}

Put $G=\mathrm{SO}_0(n+1,1)$, the connected component of the identity in
$\mathrm{SO}(n+1,1)$, $\mathfrak{g}=\mathrm{so}(n+1,1)$ and $K=\mathrm{SO}(n+1)$; let
${\mathcal H}(K)$ be the
space of $K$-finite vectors in ${\mathcal F}_\lambda(\mathbb{S}^n)$.
The main result of this note is a classification of simple and unitary
$(\mathfrak{g},K)$-submodules of~${\mathcal H}(K)$ as a function of
$\lambda$.

\begin{theorem}
\begin{enumerate}
\item
If $\lambda\not=l/n$ for $l\in\mathbb{Z}$, or if, for $n>1$,
$\lambda\in\left\{\frac{1}{n},\frac{2}{n},\dots,\frac{n-1}{n}\right\}$, then
${\mathcal F}_\lambda(\mathbb{S}^n)$ contains a unique simple $(\mathfrak{g},K)$-module
${\mathcal H}(K)$,  identified to the space of harmonic polynomials on
$\mathbb{S}^n$. This module is unitary if and only if
$\lambda=\frac{1}{2}+i\alpha$, $\alpha\in{\mathbb R}^*$, or
$\lambda\in ]0,1[\setminus\left\{\frac{1}{2}\right\}$.
\item
If $\lambda=-l/n$, $l\in\mathbb{N}$, ${\mathcal H}(K)$ contains a unique simple
$(\mathfrak{g},K)$-submodule, which is finite-dimensional and
given by the elements of degree $\le l$. It is unitary if and only
if $\lambda=0$.
\item
If $n=1$ and $\lambda=l$, $l\in\mathbb{N}^*$, ${\mathcal H}(K)$ contains two
simple $(\mathfrak{g},K)$-submodules, unitary and infinite-dimensional, and
the direct sum of these modules consists of the elements of ${\mathcal H}(K)$
of degree $\ge l$.
\item
If $n>1$ and $\lambda=1+l/n$, $l\in\mathbb{N}$, ${\mathcal H}(K)$ contains a
simple infinite-dimensional $(\mathfrak{g},K)$-submodule consisting of the
elements with degree $\ge l+1$. It is unitary if and only if $\lambda=1$.
\end{enumerate}
\label{class}
\end{theorem}

\begin{remark}
We described all the closed $G$-submodules of
${\mathcal F}_\lambda(\mathbb{S}^n)$ (cf. \cite{Kna}, Theorem~8.9), and, since $G$ is
connected,  we obtained, in the case (2), every simple finite-dimensional
$\mathfrak{g}$-submodules  of ${\mathcal F}_\lambda(\mathbb{S}^n)$.
\end{remark}

\section{Nonunitary spherical series}

The main ingredient of the proof of Theorem~\ref{class} is
the identification of the modules
${\mathcal F}_\lambda(\mathbb{S}^n)$ with induced representations. Denote $G=KAN$ the
Iwasawa decomposition of $G$ and $\rho$ the half-sum of the positive
restricted roots of the pair ($\mathrm{so}(n+1,1),\mathfrak{a})$, $A=\exp
\mathfrak{a}$.

Consider the representation ${\rm Ind}_{MAN}^G(0\otimes \nu)$, induced from
the minimal parabolic subgroup
$MAN$ of $G$, with the trivial representation of the subgroup $M=\mathrm{SO}(n)$
(the centralizer of $A$ in $K$) and a one-dimensional representation
 $\mu$ of $A$ such that, for $h\in A$, one has
$\mu(h)=\exp(\nu(\log h))$, with a fixed $\nu\in \mathfrak{a}^*$.
Abusing the notations, we identify an element $\nu$ in
$\mathfrak{a}^*$ with $\nu(H)$, where $H$ is the matricial element 
\[
H={\rm E}_{n+1,n+2}+{\rm E}_{n+2,n+1}\quad ({\rm with~elementary~matrices~} {\rm E}_{ij})
.
\]
Therefore, $\rho=\frac{n}{2}$.

The Iwasawa
decomposition shows that this induced representation
acts on the space of functions in
${\mathcal L}^2(K/M)={\mathcal L}^2(\mathbb{S}^n)$,
and the operators of this representation are given, for $g\in G$, by
\[
{\rm Ind}_{MAN}^G(0\otimes \nu)(g)f(k)=
\exp(-\nu(\log h))f(k_g),\qquad
{\rm with} \quad g^{-1}k=k_ghn\in KAN.
\]
Considering every value of $\nu$ in ${\mathbb C}$, we obtain the
representations of the so-called {\it nonunitary spherical series}, that
defines a structure of $G$-module on the space ${\mathcal L}^2(\mathbb{S}^n)$. We
denote by ${\mathcal C}^\infty_\nu(\mathbb{S}^n)$ the submodule constituted of
${\mathcal C}^\infty$ elements.

Our proof is based on the following fact.

\begin{theorem}
The $\mathfrak{g}$-modules
${\mathcal F}_\lambda(\mathbb{S}^n)$ and ${\mathcal C}^\infty_\nu(\mathbb{S}^n)$ are
isomorphic if and only if $\nu=n\lambda$, and this
isomorphism is compatible with the action of $K$.
\label{equiv}
\end{theorem}

Let us give the main idea of the proof of Theorem~\ref{equiv}. Denote by
$d\,{\rm Ind}_{MAN}^G(0\otimes \nu)$ the infinitesimal representation
associated with  ${\rm Ind}_{MAN}^G(0\otimes \nu)$,
 $L_X$ the Lie derivative along $X\in\mathfrak{g}$, and
$(\theta_1,\dots,\theta_n)$ the spherical coordinates on $\mathbb{S}^n$.
Straightforward but complicated computations lead to the following two
facts, that use cohomological (elementary) notions.

\begin{lemma}
For all $X\in\mathfrak{g}$ one has
\[
d\,{\rm Ind}_{MAN}^G(0\otimes \nu)(X)=L_X+\nu \,c(X),
\]
where $c$ is the
1-cocycle on $\mathrm{so}(n+1,1)$ with coefficients in ${\mathcal C}^\infty({\mathbb R}^n)$
given, in spherical coordinates, by
$c(X)=\frac{\partial X^n}{\partial\theta_n}$.
\end{lemma}

It is known that the cohomology space $H^1(\mathrm{so}(n+1,1);{\mathcal
C}^\infty({\mathbb R}^n))$ is one-dimensional. We then have the following

\begin{lemma}
The cocycle $c$ is cohomological to the cocycle $\tilde{c}$ given in
spherical coordinates by
\[
\tilde{c}(X)=\frac{1}{n}\,{\mathrm{Div}}\, X.
\]
\end{lemma}

We now use the fact that two representations that are given by $L_X+c(X)$
and $L_X+\tilde{c}(X)$ are equivalent if the cocycles $c$ and $\tilde{c}$
belong to the same cohomology class. Theorem~\ref{equiv} is proved.

As a consequence, $(\mathfrak{g},K)$-modules of $K$-finite vectors in
${\mathcal F}_\lambda(\mathbb{S}^n)$ and ${\mathcal C}^\infty_\nu(\mathbb{S}^n)$ are
isomorphic.

\section{Classification of $\boldsymbol{(\mathfrak{g},K)}$-modules in
$\boldsymbol{{\mathcal F}_\lambda(\mathbb{S}^n)}$}

Let us now use the results (and the notations) of \cite{Gui} (see Appendix
B.10). Let us put
$k=\left[\frac{n+1}{2}\right]$ (where $[p]$ is the integral part of $p$), and denote by
$D^{m_1,\dots,m_k}$ the simple $K$-module with highest weight
$m_{k}\varepsilon_1+m_{k-1}\varepsilon_2+\cdots+m_1\varepsilon_{k}$, where
\[
\varepsilon_i(\lambda_1 H_1+\cdots +\lambda_k H_k)=\lambda_i,
\]
and the matricial $H_r=i({\rm E}_{2r-1,2r}-{\rm E}_{2r,2r-1})$ generate a Cartan
subalgebra of the Lie algebra~$\mathfrak{k}$.

Consider the representation ${\rm Ind}_{MAN}^G(0\otimes \nu+\rho)$
(which is unitary if and only if $\nu$ is pure imaginary), and describe
the $(\mathfrak{g},K)$-module $E_{0,\nu}$ of its $K$-finite vectors~: the
restriction of the latter to $K$ is given by the direct sum of simple
$K$-modules
\[
E_{0,\nu}|_K=
\bigoplus D^{0,\dots,0,m},\,\left|\begin{array}{l}m\in\mathbb{N}\quad {\rm
for} \ \ n>1;\\m\in\mathbb{Z}\quad {\rm for} \ \ n=1.\end{array}\right.
\]
We use the isomorphism $E_{0,-\nu}\cong E_{0,\nu}^*$ ($K$-finite dual).

The module $E_{0,\nu}$ is unitary if and only if $\nu$ is pure imaginary, or
$\nu\in \left]-\frac{n}{2},\frac{n}{2}\right[\setminus\{0\}$.

In order to study simple $(\mathfrak{g},K)$-submodules of $E_{0,\nu}$, we
have to consider the following two cases.
\begin{itemize}
\item
If $n=1$, then the module $E_{0,\nu}$ is simple if and only if
$\nu\not\in\frac{1}{2}+\mathbb{Z}$.

Otherwise, we have:
\begin{itemize}
\item If $\nu<0$, $E_{0,\nu}$ contains a unique
simple$(\mathfrak{g},K)$-submodule. It is finite-dimensional and given, as a
$K$-module, by $\bigoplus_{|m|\le |\nu|-\frac{1}{2}}  D^{m}$. This module is
unitary for $\nu=-\frac{1}{2}$.

\item If $\nu>0$, $E_{0,\nu}$ contains two simple infinite-dimensional
$(\mathfrak{g},K)$-submodules, given as  $K$-modules by
$\bigoplus_{\nu+\frac{1}{2}\le \pm m}  D^{m}$. These modules are unitary.
\end{itemize}

\item
If $n>1$, then the module $E_{0,\nu}$ contains a simple submodule if and
only if
$\nu=\pm\left(\frac{n}{2},\frac{n}{2}+1,\dots\right)$. In this case, there exists a
simple finite-dimensional $(\mathfrak{g},K)$-module given, as a
$K$-module, by $\bigoplus_{m\le |\nu|-\frac{n}{2}} D^{0,\dots,0,m}$. This is
a $(\mathfrak{g},K)$-submodule of $E_{0,\nu}$ if $\nu<0$,  and a
quotient-module if
$\nu>0$. It is unitary for $\nu=-\frac{n}{2}$.
\end{itemize}

Consider the space of smooth functions on ${\mathbb R}^{n+1}\setminus
\{0\}$, homogeneous of degree $-\lambda(n+1)$.
This space is a $\mathrm{Diff}(\mathbb{S}^n)$-module (and also a
$\mathrm{Vect}(\mathbb{S}^n)$-module with the Lie derivative) isomorphic to the
module ${\mathcal F}_\lambda(\mathbb{S}^n)$ (see~\cite{Ovs}).
Denote by ${\mathcal H}^{n+1,m}$ the
$K$-module constituted of its elements of the form
\[
\frac{P_m(x_0,\dots,x_n)}
{\left(x_0^2+\dots+x_n^2\right)^{\frac{m}{2}+\lambda\frac{n+1}{2}}},
\]
where $P_m$ is a harmonic polynomial homogeneous of degree~$m$. We,
finally, check the following facts:
\begin{itemize}
\item
If $n=1$, then ${\mathcal H}^{2,m}\cong D^{-m}\oplus D^m$. Indeed ${\mathcal
H}^{2,m}$ is the direct sum of $\mathrm{SO}(2)$-modules $H_m$ and $H_{-m}$,
respectively generated by
\[
(x_0+ix_1)^m\left(x_0^2+x_1^2\right)^{-\frac{m}{2}-\lambda}
\]
and
its conjugate in ${\mathbb C}$, and we have $H_{\pm m}\cong D^{\pm m}$.

\item
If $n>1$, then ${\mathcal H}^{n+1,m}$ is simple and we have ${\mathcal
H}^{n+1,m}\cong D^{0,\dots,0,m}$ .
\end{itemize}

Consequently, the $(\mathfrak{g},K)$-module of $K$-finite
vectors of ${\mathcal F}_\lambda(\mathbb{S}^n)$ is given by
\[
{\mathcal H}(K)\cong\bigoplus_{m\in\mathbb{N}}{\mathcal H}^{n+1,m}.
\]

Let us apply the above results to the representation ${\rm
Ind}_{MAN}^G(0\otimes \nu)$. Substituting $\nu+\rho=n\lambda$ to
Theorem~\ref{equiv}, we obtain the assertions of Theorem~\ref{class}.

\begin{remark}
The case $n=1$ can be directly deduced from the classification of
representations of $\mathrm{SL}(2,{\mathbb R})$: acting the same way as in~\cite{Lan}, we
observe that the space of $K$-finite vectors of
${\mathcal F}_\lambda(\mathbb{S}^1$) is the direct sum $\bigoplus
H_{2l},~l\in\mathbb{Z}$, where $H_m$  is the space of the representation of
$\mathrm{SO}(2)\subset\mathrm{SL}(2,{\mathbb R})$ with the character
\[
\chi_m:
\left[\begin{array}{cc} \cos\theta&\sin\theta\\
-\sin\theta&\cos\theta\end{array}
\right]\mapsto{\rm e}^{im\theta}.
\]
\end{remark}

\subsection*{Acknowledgments}

We are grateful to Valentin Ovsienko, Thierry Levasseur and Alain
Guichardet, so as to Ranee Brylinski, Patrick Delorme and Pierre Lecomte.

\label{redou-lastpage}

\end{document}